\documentclass{amsart}

\usepackage{amsmath,amssymb,amsthm}
\usepackage{graphicx,tikz}

\newtheorem*{thm}{Theorem}

\theoremstyle{definition}

\theoremstyle{remark}

\begin{document}

\title[]{A forgotten Theorem of Sch\"onberg\\ on one-sided integral averages}

\keywords{Integral averages, aggregate function, exponential distribution.}
\subjclass[2010]{44A35, 62M10}

\author[]{Stefan Steinerberger}
\address{Department of Mathematics, Yale University, New Haven, CT 06511, USA}
\email{stefan.steinerberger@yale.edu}
\thanks{S.S. is supported by the NSF (DMS-1763179) and the Alfred P. Sloan Foundation.}

\begin{abstract}  Let $f:\mathbb{R} \rightarrow \mathbb{R}$ be a function for
which we want to take local averages. Assuming we cannot look into the future, the 'average'
at time $t$ can only use $f(s)$ for $s \leq t$.
 A natural way to do so is via a weight $\phi$ and 
$$ g(t) = \int_{0}^{\infty}{f(t-s) \phi(s) ds}.$$
We would like that (1) constant functions, $f(t) \equiv \mbox{const}$, are mapped to themselves and (2) $\phi$ to be monotonically decreasing (the more recent past should weigh more heavily than the
distant past). Moreover, we want that (3) if $f(t)$ crosses a certain threshold $n$ times, then $g(t)$ should not cross
the same threshold more than $n$ times (if $f(t)$ is the outside wind speed and crosses the Tornado threshold at two points in time, we would like the averaged wind speed to cross the Tornado threshold \textit{at most} twice). A Theorem implicit in the work of Sch\"onberg is that these three conditions characterize a unique weight that is given by the exponential distribution
$$ \phi(s) = \lambda^{} e^{-\lambda s} \qquad \mbox{for some} \quad \lambda > 0.$$
\end{abstract}

\maketitle

\section{Introduction and Result}
The purpose of this paper is to discuss how one would go about averaging continuous functions.
 Let $f:\mathbb{R} \rightarrow \mathbb{R}$ be a continuous function and suppose
that we are interested in, at a given time $t$, finding a local average of $f$ using only function values $f(s)$ for $s \leq t$.
This is the canonical setting for many applications where we cannot look into the future (one only needs to think of sports or finance where this is a constant problem).
A natural way of constructing an average is via 
$$ g(t) = \int_{0}^{\infty}{f(t-s) \phi(s) ds},$$
where $\phi:[0,\infty] \rightarrow \mathbb{R}$ is a (not necessarily continuous) weighting function. Many different weighting functions are conceivable, the one that is presumably used most often in practice is
$$ \phi(s) = \begin{cases} a^{-1} \qquad &\mbox{if}~0 \leq s \leq a \\ 0 \qquad &\mbox{otherwise,} \end{cases}$$
the average taken over the last $a$ units of time. A natural question is whether there is a 'best' weight and, as usual, this depends on how one defines things. We will proceed in an axiomatic fashion and state a list of desirable properties.

\begin{quote}
\textbf{Property 1. Invariance of Constants.} Averaging should leave constant functions, $f(t) \equiv c$, invariant.
\end{quote}
\begin{quote}
\textbf{Property 2. Monotonicity.} $\phi:[0,\infty] \rightarrow \mathbb{R}$ is (not necessarily strictly) monotonically decreasing.
\end{quote}
\begin{quote}
\textbf{Property 3. Variation-diminishing property.} For any $c \in \mathbb{R}$, if the set $\left\{t \in \mathbb{R}: f(t) > c\right\}$ is a union of $n$ (not necessarily bounded) intervals, then the set $\left\{t \in \mathbb{R}: g(t) > c\right\}$ is the union of at most $n$ intervals. If $\left\{t \in \mathbb{R}: f(t) < c\right\}$ is the union of at most $n$ (not necessarily bounded) intervals, then so is  $\left\{t \in \mathbb{R}: g(t) < c\right\}$. \end{quote}

The first condition is completely unambiguous: an average of constant values needs to return the same value in order to be meaningful. This translates easily into
$$ \int_{0}^{\infty}{\phi(s) ds} = 1.$$
We observe that this condition also implies that any function $\phi$ satisfying Properties (1) and (2) is nonnegative: if it assumes negative values anywhere, then monotonicity would imply that it is not integrable which violates Property (1). In particular, $\phi$ is a probability distribution.
As a consequence of that, we have that
$$ \min_{s \leq t}{ f(s) } \leq g(t) = \int_{0}^{\infty}{f(t-s) \phi(s) ds} \leq \max_{s \leq t}{f(s)}$$
which is also exceedingly natural: the average value at any point cannot exceed the previously attained maximal value or be smaller than all previous values.
The second condition, monotonicity, is natural insofar as we would like the recent past to be more representative than the distant past.
Property (3) is a smoothing property: the averaged function should not venture into 'extreme' territory  $\left\{t \in \mathbb{R}: g(t) > c\right\}$ more often than the function does itself $\left\{t \in \mathbb{R}: f(t) > c\right\}$. Requiring property (3) to be satisfied for all $c \in \mathbb{R}$ therefore corresponding to a uniform smoothing at all scales -- \textit{extreme events can be represented in the average but they should not be over-represented.}
A simple example is as follows: suppose $f(t)$ is the ELO strength of a chess player measured at time $t$. This indicator is discontinuous and changes after each game -- however, if a chess players has their ELO exceed 2800 for the entirety of the year 2016 and then once more, briefly, in 2018, then it would be desirable for the averaged function to exceed the value 2800 at most two times and not, say, three times. It would be perfectly reasonable, however, if the averaged function $g$ exceeds the value 2800 only once or never at all (for example if the value in 2016 hovers very close to 2800 all the time and was much lower before or, conversely, if in the month of 2018 the value is only exceeded for a brief period of time). \\

These three conditions uniquely characterize a weight (up to dilation symmetries).

\begin{thm}[Sch\"onberg] If a function $\phi$ satisfies properties (1), (2) and (3), then
$$ \phi(s) = \lambda ^{} e^{-\lambda s} \qquad \mbox{for some}~\lambda > 0.$$
\end{thm}
To the best of our knowledge, this Theorem has never been stated or proved. I. J. Sch\"onberg mentions in passing in his 1948 paper \cite{sch2}
that, as a consequence of his classification theorem, 'All Polya frequency functions turn out to be continuous everywhere
with the single exception of the truncated exponential' and this is exactly what is needed to prove the Theorem which should be attributed to Sch\"onberg. The use of exponential distributions $\phi$ to compute one-sided averages is completely classical in time series analysis ('exponential smoothing') and usually ascribed to work of Brown \cite{brown} or Holt \cite{holt} in the 1950s but the fact that properties (1) -- (3) uniquely characterize exponential smoothing does not seem to be known.\\

We emphasize that, as one often encounters in axiomatic approaches, the result is only as good as one's faith in the axioms. This is the second purpose of this paper: to perhaps motivate a study of axiomatic approaches towards integral averages. What properties should an integral averaging operator have and which types of averages possess these properties? We believe all three properties to fairly natural (with (3) being a particularly subtle way of defining smoothing). As is customary in axiomatic approaches, there are presumably other axioms that might also be of interest and will generally lead to different results. 
\section{The Proof}
\begin{proof} As discussed above, properties (1) and (2) imply that $\phi$ is a probability density function. This implies, by linearity, that the function is invariant under adding constants. This allows us to replace the study of
$$ \left\{ t \in \mathbb{R}: f(t) \geq c \right\} \qquad \mbox{and} \qquad  \left\{ t \in \mathbb{R}: g(t) \geq c \right\} $$
with the study of when $f$ and $g$ become positive (by replacing $f(t)$ with $f(t)-c$ which leads to $g(t)$ being replaced by $g(t) - c$).
Property (3) is then equivalent to asking that the number of sign changes of $g$ is at most that of the number of sign changes of $f$ or equivalently, it asks that convolution with $\phi$ has the variation-diminishing property in the sense of
Sch\"onberg. Phrased differently, we learn that $\phi$ is a Polya frequency function \cite{deb, spl1, mar, pol1, pol2}.
Sch\"onberg's theory \cite{sch, sch2, sch3, sch4} implies 
$$ \int_{\mathbb{R}}{ e^{-sx} \phi(x) dx} = \frac{1}{\psi(s)} \qquad \mbox{for all complex}~s~\mbox{in}~ -a < \Re s < a$$
for some $a > 0$ where the function $\psi(s)$ is an entire function of the form
$$ \psi(s) = C e^{-\gamma s^2 + \delta s} \prod_{k=1}^{\infty}{(1+\delta_k s) e^{-\delta_k s}},$$
where $C>0$, $\gamma \geq 0$, $\delta$ and $\delta_k$ are real numbers and $\sum_{k=1}^{\infty}{\delta_k^2} < \infty$. It
remains to find all functions $\phi$ of that type satisfying all our properties. We analyze the behavior for purely imaginary $s = it$ where $t \in \mathbb{R}$.
Since $\gamma, \delta$ and $\delta_k$ are real,
$$ | \psi(it)| = C e^{\gamma t^2}\prod_{k=1}^{\infty}{|1 + \delta_k i t|} \geq C\prod_{k=1}^{\infty}{|1 + \delta_k i t|}.$$
The product is well defined since, using $\log{(1+x)} \leq x$ for $x>0$,
$$ \log \left( \prod_{k=1}^{\infty}{|1 + \delta_k i t|} \right) = \sum_{k=1}^{\infty}{ \log{\left(\sqrt{1+ \delta_k^2 t^2}\right)}} \leq \frac{t^2}{2} \sum_{k=1}^{\infty}{\delta_k^2}.$$
Moreover, we have $\|1 + \delta_k t\| \geq 1$ for all $k \in \mathbb{N}$. We distinguish two cases: either all but one $\delta_k$ are zero or at least two $\delta_k$ are nonzero. In the second case, we observe  
$$ | \psi(it) | \geq 1 + c t^2$$
for some fixed $c>0$. Applying the inverse Fourier transform shows that
$$ \phi(s) = \int_{\mathbb{R}}{ e^{isx} \frac{1}{\psi(x)}dx}.$$
We note that $\phi(0) > 0$ by assumption (2) as well as $\phi(s) = 0$ for all $s < 0$ by assumption.
Let $R \gg 1$ be so large that
$$ \int_{\mathbb{R} \setminus [-R,R]}{ \frac{1}{|\psi(x)|}dx} \leq \frac{\phi(0)}{4}.$$
We then observe that
\begin{align*}
\phi(0) &= \left| \phi(-R^{-2}) - \phi(0)\right| \\
&= \left| \int_{\mathbb{R}}{ (e^{-R^{-2}ix} - 1) \frac{1}{\psi(x)}dx}\right|\\
&\leq \left| \int_{-R}^{R}{ (e^{-R^{-2}ix} - 1) \frac{1}{\psi(x)}dx}\right| + \left| \int_{\mathbb{R} \setminus [-R,R]}^{}{ (e^{-R^{-2}ix} - 1) \frac{1}{\psi(x)}dx}\right| \\
&\leq \frac{1}{R} \int_{\mathbb{R}}{ \frac{1}{|\psi(x)|}dx} + \frac{\phi(0)}{2}.
\end{align*}
The first term can be made sufficiently small by further increasing $R$. This contradiction shows that 
all but exactly one $\delta_k$ to be 0. The same argument can be run of $\gamma > 0$. Thus
$$ \psi(s) = C e^{\delta s} (1+ \delta s) e^{-\delta s} = C (1+\delta s).$$
This shows that
$$ \int_{0}^{\infty}{e^{-i s x} \phi(x) dx} = \frac{C}{1+ \delta i s} \qquad \mbox{for all}~s \geq 0.$$
We define a function 
$$ \phi_2(x) = \begin{cases} \phi(x) \qquad &\mbox{if}~x \geq 0 \\ \phi(-x) \qquad &\mbox{if}~x \leq 0. \end{cases}$$
We note that
\begin{align*} \int_{\mathbb{R}}^{}{e^{-i s x} \phi_2(x) dx} &= \int_{0}^{\infty}{(e^{-i s x}+e^{i s x}) \phi(x) dx} \\
&= \int_{0}^{\infty}{2 \cos{(sx)} \phi(x) dx} = 2 \Re \frac{C}{1+ \delta i s} = \frac{2 C}{1+ \delta^2 s^2}.
\end{align*}
An application of the inverse Fourier transform shows that
$$ \phi_2(x) = C e^{-\frac{|x|}{\delta}}.$$
The normalization
$$ \int_{0}^{\infty}{\phi(x) dx} = 1$$
then implies
$$ \phi(x) = \lambda^{} e^{-\lambda x} \qquad \mbox{for some}~\lambda > 0.$$
\end{proof}

\end{document}